\documentclass[11pt]{article}
\usepackage{amsfonts, amsmath, amssymb, amsthm}

\usepackage[colorlinks=true,citecolor=blue,urlcolor=blue,linkcolor=blue,bookmarksopen=true]{hyperref}

\title{Bounds on equiangular lines and on related spherical codes}
\author{Boris Bukh\footnote{Supported in part by U.S.\ taxpayers through NSF grant DMS-1301548. Supported
in part by Alfred P.\ Sloan Foundation through Sloan Research Fellowship.}}
\date{}

\makeatletter
\def\TODO{\@ifnextchar[{\TODO@with}{\marginnote{TODO}}}
\def\TODO@with[#1]{\marginnote{#1}}
\makeatother

\newtheorem{theorem}{Theorem}
\newtheorem{lemma}[theorem]{Lemma}

\newtheorem{conjecture}[theorem]{Conjecture}

\newcommand*{\abs}[1]{\lvert #1\rvert}                           % Absolute values, cardinality
\newcommand*{\norm}[1]{\lVert #1\,\lVert}                        % Norm, l^1 norm
\newcommand*{\eqdef}{\stackrel{\text{\tiny{def}}}{=}}            % definition by equality                                      
\newcommand*{\veps}{\varepsilon}                                 % Nicely-looking epsilon
\newcommand*{\R}{\mathbb{R}}                                     % Complex numbers
\DeclareMathOperator{\vspan}{span}                               % Span of a set of vectors

\begin{document}

\maketitle

\begin{abstract}
An $L$-spherical code is a set of Euclidean unit vectors whose pairwise
inner products belong to the set~$L$. We show, for a fixed $0<\alpha,\beta<1$, 
that the size of any $[-1,-\beta]\cup\{\alpha\}$-spherical code is at most
linear in the dimension.

In particular, this bound applies to sets of lines such that every two
are at a fixed angle to each another.
\end{abstract}

\section{Introduction}
\paragraph{Background}
A set of lines in $\R^d$ is called \emph{equiangular}, if the angle between any two of them
is the same. Equivalently, if $P$ is the set of unit direction vectors, the
corresponding lines are equiangular with the angle $\arccos \alpha$ if $\langle v,v'\rangle\in\{-\alpha,\alpha\}$ for any
two distinct vectors $v,v'\in P$. The second equivalent way of defining equiangular
lines is via the Gram matrix. Let $M$ be the matrix whose columns are the direction vectors. Then
$M^TM$ is a positive semidefinite matrix whose diagonal entries are $1$'s, and each of whose off-diagonal entries is $-\alpha$
or~$\alpha$. Conversely, any such matrix of size $m$ and rank $d$ gives rise to $m$ equiangular lines in~$\R^d$.

Equiangular lines have been extensively studied following the works of van Lint and Seidel \cite{van_lint_seidel}, and of Lemmens and Seidel \cite{lemmens_seidel}.
Let $N(d)$ be the maximum number of equiangular lines in $\R^d$. Let $N_{\alpha}(d)$ be the maximum
number of equiangular lines with the angle $\arccos \alpha$. The values of $N(d)$ are known exactly for $d\leq 13$, for $d=15$, for $21\leq d\leq 41$ and for $d=43$ \cite{barg_yu,greaves_koolen_munemasa_szollosi}. When $d$ is large, the only known upper bound on $N(d)$ is due to
Gerzon (see \cite[Theorem 3.5]{lemmens_seidel}) and asserts that
\[
  N(d)\leq d(d+1)/2\ \text{with equality only if }d=2,3\text{ or }d+2\text{ is a square of an odd integer}.
\]
A remarkable construction of de Caen\cite{decaen_equiangular} shows that $N(d)\geq \tfrac{2}{9}(d+1)^2$ for $d$ of the
form $d=6\cdot 4^i-1$. A version of de Caen's construction suitable for other values of $d$ 
has been given by Greaves, Koolen, Munemasa and Sz\"oll\"osi \cite{greaves_koolen_munemasa_szollosi}.
See also the work of Jedwab and Wiebe \cite{jedwab_wiebe_equiangular} for an alternative construction of $\Theta(d^2)$ equiangular lines.
In these constructions the inner product $\alpha$ tends to $0$ as dimension grows.

\paragraph{Previously known bounds on $N_{\alpha}(d)$}
The first bound is the so-called \emph{relative bound} (see \cite[Lemma~6.1]{van_lint_seidel} following \cite[Theorem~3.6]{lemmens_seidel})
\[
  N_{\alpha}(d)\leq d\frac{1-\alpha^2}{1-d\alpha^2}\quad\text{if }d<1/\alpha^2.
\]
While useful in small dimensions, it gives no information for a fixed $\alpha$ and large~$d$. The second bound
is
\[
  N_{\alpha}(d)\leq 2d\quad\text{unless }1/\alpha\text{ is an odd integer \cite[Theorem 3.4]{lemmens_seidel}}.
\]
This bound can be further improved to $\tfrac{3}{2}(d+1)$ unless $\frac{1}{2\alpha}+\frac{1}{2}$ is an algebraic integer of degree~$2$,
see \cite[Subsection 2.3]{bukh_lmatrices}. 

Finally, the values of $N_{1/3}(d)$ and $N_{1/5}(d)$ for a large $d$ have been completely determined:
\begin{align*}
  N_{1/3}(d)&=2d-2&&\text{for }d\geq 15&\quad&\text{\cite[Theorem~4.5]{lemmens_seidel}},\\
  N_{1/5}(d)&=\lfloor 3(d-1)/2\rfloor&&\text{for all sufficiently large }d&\quad&\text{\cite{neumaier} and \cite[Corollary 6.6]{greaves_koolen_munemasa_szollosi}}.
\end{align*}

\paragraph{New bound} We will show that $N_{\alpha}(d)$ is linear for every $\alpha$. In fact, we will
prove a result in greater generality. Following \cite{delsarte_goethals_seidel}, we call a set of unit 
vectors $P$ an \emph{$L$-spherical code} if $\langle v,v'\rangle \in L$ for every pair of distinct points 
$v,v'\in P$. In this language, a set of equiangular lines is a $\{-\alpha,\alpha\}$-spherical code.
Let $N_L(d)$ be the maximum cardinality of an $L$-spherical code in~$\R^d$.

\begin{theorem}\label{thm:main}
For every fixed $0<\beta\leq 1$ there exists a constant $c_{\beta}$ such that for any $L$ of the form $L=[-1,-\beta]\cup\{\alpha\}$
 we have $N_L(d)\leq c_{\beta} d$.
\end{theorem}

We make no effort to optimize the constant $c_{\beta}$ that arises from our proof, as it is huge. We speculate
about the optimal bounds on $N_L(d)$ in section~\ref{sec:open}. We do not know if the constant $c_{\beta}$ 
can be replaced by an absolute constant that is independent of $\beta$, i.e., whether 
$N_L(d)\leq cd+o_{\beta}(d)$ holds.

The rest of the paper is organized as follows. In the next section we prove Theorem~\ref{thm:main} and
in the concluding section we discuss possible generalizations and strengthenings of Theorem~\ref{thm:main}.

\section{Proof of Theorem~\protect\ref{thm:main}}
\paragraph{Proof sketch}
The idea behind the proof of Theorem~\ref{thm:main} builds upon the argument of Lemmens and Seidel for $N_{1/3}(d)$.
Before going into the details, we outline the argument. 

Let $L=[-1,-\beta]\cup\{\alpha\}$, and let $P$ be an $L$-spherical code whose size we wish to bound.
Define a graph $G$ on the vertex set $P$ by connecting $v$ and $v'$ by an edge if 
$\langle v,v'\rangle\in [-1,-\beta]$. In their treatment of $N_{1/3}(d)$ Lemmens and Seidel consider
the largest clique in $G$, and carefully analyze how the rest of the graph attaches to that clique.
In contrast, in our argument we consider the largest independent set $I$ of $G$, and show that almost
every other vertex is incident to nearly all vertices of~$I$. Iterating this argument
inside the common neighborhood of $I$ we can build a large clique in $G$. As the clique size is bounded
by a function of $\beta$, that establishes the theorem.

\paragraph{Proof details} For the remainder of the section, $L$, $P$ and $G$ will be as defined as in the preceding proof sketch.
The following two well-known lemmas bound the sizes of cliques and independent sets in $G$:

\begin{lemma}\label{lem:negative}
Suppose $u_1,\dotsc,u_n$ are $n$ vectors of norm at most $1$ satisfying $\langle u_i,u_j\rangle\leq -\gamma$.
Then $n\leq 1/\gamma+1$.
\end{lemma}
\begin{proof}
This follows from $0\leq \norm{\sum u_i}^2=\sum_{i,j}\langle u_i,u_j\rangle
\leq n-\gamma n(n-1)$.
\end{proof}
\begin{lemma}\label{lem:homog}\
\begin{enumerate}
\item Every independent set in $G$ is linearly independent. In particular, the graph $G$ contains no independent set on more than $d$ vertices.
\item The graph $G$ contains no clique on more than $1/\beta+1$ vertices.
\end{enumerate}
\end{lemma}
\begin{proof}
i) 
Let $p_1,\dotsc,p_n$ be the points of the independent set. Suppose that $\sum c_i p_i=0$. Taking an inner product
with $p_j$ we obtain $0=(1-\alpha)c_j+\alpha\sum c_i$ implying that all $c_i$'s are equal. The result follows since
$(1-\alpha)+n\alpha\neq 0$.

ii) This is a special case of the preceding lemma.
\end{proof}

In the next two lemmas we analyze how the vertices of $G$ attach to an independent set.
\begin{lemma}\label{lem:dist}
Suppose that $M$ is a matrix with linearly independent column vectors $p_1,\dotsc,p_n$.
Suppose that $v,v'\in \vspan \{p_1,\dotsc,p_n\}$ are points satisfying $\langle p_i,v\rangle=s_i$
and $\langle p_i,v'\rangle=s_i'$ for some column vectors $s=(s_1,\dotsc,s_n)^T$ and $s'=(s_1',\dotsc,s_n')^T$. 
Then $\langle v,v'\rangle=s^T(M^TM)^{-1}s'$.
\end{lemma}
\begin{proof}
By passing to a subspace we may assume that $p_1,\dotsc,p_n$ span $\R^n$, and so
$M$ is invertible. As $s=M^Tv$ and $s'=M^Tv'$, we infer that $\langle v,v'\rangle=v^Tv'=((M^T)^{-1}s)^T(M^T)^{-1}s'=s^T(M^TM)^{-1}s'$.
\end{proof}

The following lemma is the geometric heart of the proof. In its special case $v=v'$, the lemma bounds degrees
from certain vertices into an independent set. More precisely, let $I$ be a sufficiently large independent set.
We will show later (in Lemma~\ref{lem:bad}) that the vertices, the norm of whose projection on $\vspan I$ exceeds $\alpha^{1/2}$, 
are few. The straightforward, but slightly messy calculations in the following lemma
characterize the vertices with such projections in terms of their degree into~$I$. 
The case $v\neq v'$ is not needed when $P$ comes from a set of equiangular lines,
but is required to establish Theorem~\ref{thm:main} in its full generality.
%The extra generality, that is gained by allowing $v\neq v'$, is necessary to handle the 
%case when $P$ does not come from a set of equiangular lines, but from the more general
%$L$-spherical codes of Theorem~\ref{thm:main}.
\begin{lemma}\label{lem:geo}
Let $t=1/\beta+1$.
There exists $n_0=n_0(\beta)$ and
$\veps=\veps(\beta)>0$ such that the following holds.
Suppose that $p_1,\dotsc,p_n$ is an independent set in $G$ of size $n$, and suppose that points $p,p'\in P$ are adjacent to
the same $m$ vertices among $p_1,\dotsc,p_n$. Assume $0<m<n-t$ and $n\geq n_0$.
Write $p=v+u$ and $p'=v'+u'$ where $v,v'\in \vspan \{p_1,\dotsc,p_n\}$ and $u,u'$ are both orthogonal to $\vspan \{p_1,\dotsc,p_n\}$.
Then $\langle v,v'\rangle\geq \alpha+\veps$.
\end{lemma}
\begin{proof}
For the duration of this proof, $I$ denotes the identity matrix, and $J$ denotes the all-$1$ matrix.
Let $M$ be the matrix comprised of column vectors $p_1,p_2,\dotsc,p_n$.
Since points $p_1,\dotsc,p_n$ are linearly independent (by Lemma~\ref{lem:homog}), the condition of the preceding
lemma is fulfilled. We have $M^TM=\alpha J+(1-\alpha)I$. One can verify that its inverse is given by
\begin{align}\label{eq:inverse}
(1-\alpha)(M^TM)^{-1}=I-\phi J\qquad\text{ with }\qquad\phi\eqdef\frac{\alpha}{1+(n-1)\alpha}.
\end{align}
%% THE FOLLOWING STEP IS UNNECESSARY. WEAKER \phi<=1/n SUFFICES
%Note that since $G$ contains a non-trivial independent set, we must have $\alpha<1$, and 
%hence $\phi<1/n$.
Note that $\phi\leq 1/n$ since $\alpha\leq 1$.

Without loss of generality, $p_1,\dotsc, p_m$ are the $m$ vertices that $p$ and $p'$ are adjacent to.
This means that $s\eqdef M^Tv$ and $s'\eqdef M^T v'$ are of the form 
$s=(-\beta_1,\dotsc,-\beta_m,\alpha,\dotsc,\alpha)$ and $s'=(-\beta_1',\dotsc,-\beta_m',\alpha,\dotsc,\alpha)$
for some $\beta_1,\beta_1',\dotsc,\beta_m,\beta_m'\in [\beta,1]$.
From Lemma~\ref{lem:dist} and \eqref{eq:inverse} it follows that
\begin{equation}\label{eq:optim}
  (1-\alpha)\langle v,v'\rangle=\alpha^2(n-m)+\sum_{i=1}^m \beta_i\beta_i'-\phi\bigl(\sum_{i=1}^n s_i\bigr)\bigl(\sum_{i=1}^n s_i'\bigr).
\end{equation}
We claim that, subject to the constraint $\beta_1,\beta_1',\dotsc,\beta_m,\beta_m'\in [\beta,1]$, the right side of \eqref{eq:optim} is minimized
when all the $\beta_i$'s and all the $\beta_i'$'s are equal to $\beta$. Indeed, since $[\beta,1]^{2m}$ is compact, the minimum is actually attained.
Assume that $(\beta_1,\beta_1',\dotsc,\beta_m,\beta_m')$ is the vector achieving the minimum, and
let $j$ be the index for which $\beta_j'$ is the largest.
Then the derivative of the right side of \eqref{eq:optim} with respect to $\beta_j'$ is 
\[
  \beta_j-\phi\sum \beta_i+(n-m)\alpha\geq \beta_j-\frac{1}{n}\sum \beta_i>\beta_j-\frac{1}{m}\sum \beta_i\geq 0.
\]
By the optimality assumption on the vector $(\beta_1,\beta_1',\dotsc,\beta_m,\beta_m')$ this implies 
that $\beta_j'=\beta$. From the choice of $j$ it then follows that $\beta_i'=\beta$ for all~$i$. Similarly, $\beta_i=\beta$ for all~$i$.

We thus deduce that
\[
  (1-\alpha)\langle v,v'\rangle\geq \alpha^2(n-m)+\beta^2m-\phi\left((n-m)\alpha-m \beta\right)^2.
\]
Let $R(m,n)$ denote the right side of preceding inequality.
Let $t^*\eqdef \frac{(1-\alpha)(\alpha-\beta)}{\alpha(\alpha+\beta)}$. We have
\[
  t^*=\frac{(1-\alpha)(\alpha-\beta)}{\alpha(\alpha+\beta)}< \frac{1-\alpha}{\alpha+\beta}< \frac{1}{\beta}=t-1.
\]
Thus to prove the lemma, it is enough to show that $R(m,n)\geq (1-\alpha)\alpha+\veps$ whenever $1\leq m\leq n-t^*-1$
and $n\geq n_0$ for suitable $n_0$ and~$\veps$. 

The expression $R(m,n)$ is a quadratic polynomial in $m$. A simple calculation shows that it satisfies $R(m,n)=R(n-t^*-m,n)$, and in particular that the maximum
of $R(m,n)$ for a fixed $n$ is at the point $m_{\max}\eqdef (n-t^*)/2$, which is inside the interval $[1,n-t^*-1]$. Furthermore, at the boundary points of the interval
we have
\begin{align*}
R(1,n)=R(n-t^*-1,n)=\alpha(1-\alpha)+(\alpha+\beta)^2-\frac{\alpha(1+\beta)^2}{1+\alpha(n-1)}.
\end{align*}
%Let $n_0=1+8/\beta^2$. Since $\frac{\alpha(1+\beta)^2}{1+\alpha(n-1)}<\frac{4}{n-1}$,
%it follows that $R(m,n)\geq R(1,n)>\alpha(1-\alpha)+\tfrac{1}{2}(\alpha+\beta)^2$ whenever $1\leq m\leq n-t^*-1$ and $n\geq n_0$.
Let $n_0=1+8/\beta^2$. Since $\frac{\alpha(1+\beta)^2}{1+\alpha(n-1)}<\frac{(1+\beta)^2}{n-1}\leq \frac{4}{n-1}$,
for $n\geq n_0$ and $1\leq m\leq n-t^*-1$ we have the inequality
$R(m,n)\geq R(1,n)>\alpha(1-\alpha)+\tfrac{1}{2}(\alpha+\beta)^2$.
In particular $\langle v,v'\rangle>\alpha+\veps$ holds under the same conditions on $n$ and $m$, where~$\veps=\tfrac{1}{2}\beta^2$.
\end{proof}

\begin{lemma}\label{lem:bad}
Suppose $p_1,\dotsc,p_n$ is an independent set in $G$. Suppose $p^{(1)},\dotsc,p^{(m)}\in P$ are points
of the form $p^{(i)}=v^{(i)}+u^{(i)}$ with $v^{(i)}\in \vspan\{p_1,\dotsc,p_n\}$ and $u^{(i)}\bot \vspan \{p_1,\dotsc,p_n\}$ and $\langle v^{(i)},v^{(j)}\rangle>\alpha+\veps$ for
all $i,j$. Then
$m\leq 1/\veps+1$.
\end{lemma}
\begin{proof}
From $\langle p^{(i)},p^{(j)}\rangle=\langle v^{(i)},v^{(j)}\rangle+\langle u^{(i)}, u^{(j)}\rangle$ and $\langle p^{(i)},p^{(j)}\rangle\in [-1,-\beta]\cup \{\alpha\}$,
we deduce that $\langle u^{(i)}, u^{(j)}\rangle<\nobreak-\veps$. The result then follows from Lemma~\ref{lem:negative}.
\end{proof}

The combinatorial part of the argument is contained in the next result.
\begin{lemma}\label{lem:ramsey}
Suppose $\delta>0$ is given. Then there exists a constant $M(\beta,\delta)$ such that the following holds.
Let $U\subset P$ be arbitrary. Suppose $I$ is a maximum-size independent subset of $U$. Then
there is a subset $U'\subset U\setminus I$ of size $\abs{U'}\geq \abs{U}-M\abs{I}$ such that 
every vertex of $U'$ is adjacent to at least $(1-\delta)\abs{I}$ vertices of~$I$.
\end{lemma}
\begin{proof}
Let $t$, $\veps$ and $n_0$ be as in Lemma~\ref{lem:geo}, and put $n=\max(n_0,\lceil t/\delta\rceil)$. Denote by $R$ the
least integer such that every graph on $R$ vertices contains either an independent set of size $n+1$ or a clique of size
at least $1/\beta+2$ (such an $R$ exists by Ramsey's theorem; furthermore, it satisfies $R\leq \binom{n+1/\beta+1}{n}$). Let 
\begin{align*}
  M&=\max(R,(1/\veps+1)2^n),\\
  N&=\abs{I}.
\end{align*}
If $\abs{U}<M$, then $\abs{U}-M\abs{I}$ is negative, and the lemma is vacuous. So, assume $\abs{U}\geq M$.
In particular, $\abs{U}\geq R$, and since by Lemma~\ref{lem:homog} the set $U$ contains no clique of size greater $1/\beta+1$,
we conclude that $N\geq n+1$.

Arrange the elements of $I$ on a circle, and consider all $N$ circular intervals containing $n$ vertices 
of $I$. Let $S_1,S_2,\dotsc,S_N$ be these intervals, in order.

We declare a vertex $p\in U\setminus I$ to be \emph{$i$-bad} if it is adjacent to between $1$ and $n-t$ vertices of $S_i$. 
For a set $T\subset S_i$, we call an $i$-bad vertex $p$ to be \emph{of type} $T$ if $T$ is precisely the set of neighbors of $p$
in the set~$S_i$. Let $B_{i,T}$ be the set of all $i$-bad vertices of type $T$, and let $B_i=\bigcup_T B_{i,T}$ be the
set of all $i$-bad vertices. By Lemmas~\ref{lem:geo} and~\ref{lem:bad} we have $\abs{B_{i,T}}\leq 1/\veps+1$ for
every $T$, and so
\[
  \abs{B_i}\leq (1/\veps+1)(2^n-1).
\]
Let $B=\bigcup B_i$ be the set of bad vertices. Hence, $\abs{B}\leq N(1/\veps+1)(2^n-1)$, and
$\abs{B\cup I}\leq MN$.\smallskip

Consider a vertex $p\in U\setminus I$ that is good, i.e., $p\not \in B$. Since $I$ is a maximal independent set,
$p$ is adjacent to at least one vertex of $I$. Say $p$ is adjacent to a vertex of $S_i$ for some $i$.
Since $p$ is good, $p$ must in fact be adjacent to at least $n-t$ vertices of $S_i$. As $S_i$
shares $n-1$ vertices with both $S_{i-1}$ and $S_{i+1}$, we are impelled to conclude that $p$ must
be adjacent to some of the vertices of $S_{i-1}$ and of $S_{i+1}$. Repeating this argument we conclude
that $p$ is non-adjacent to at most $t$ elements from among any interval of length~$n$. In particular,
$p$ is adjacent to at least $N(1-t/n)$ vertices of $I$. As $p$ is an arbitrary good vertex and $t/n\leq\delta$, the lemma follows.
\end{proof}

We are now ready to complete the proof of Theorem~\ref{thm:main}. Indeed, with foresight we set 
\begin{align*}
B&=\lceil 1/\beta+1\rceil,\\
\delta&=1/(B+1)^2.
\end{align*}
 and let $M$ be as in the proceeding lemma. Put $U_0=P$ and let $I_0$
be a maximal independent set in $U_0$. By the preceding lemma, there exists $U_1\subset U_0\setminus I_0$
such that every vertex of $U_1$ is adjacent to $(1-\delta)\abs{I_0}$ vertices of $I_0$ and $\abs{U_1}\geq \abs{U_0}-M\abs{I_0}$.
In view of Lemma~\ref{lem:homog}, $\abs{U_1}\geq \abs{U_0}-Md$. Let $I_1$ be a maximal independent
set in $U_1$. Repeating this argument, we obtain a nested sequence of
sets $U_0\supset U_1\supset $ and a corresponding sequence of independent sets $I_0,I_1,\dotsc$ such that
\begin{enumerate}
  \item $\abs{U_i}\geq \abs{U_{i-1}}-Md$ for each $i=1,2,\dotsc$,
  \item For $r<s$, each vertex in $I_s$ is adjacent to at least $(1-\delta)\abs{I_r}$ vertices of $I_r$.
\end{enumerate}

We claim that $\abs{P}\leq BMd$, which would be enough to complete the proof of Theorem~\ref{thm:main}. Indeed, suppose for the sake
of contradiction that $\abs{P}>BMd$. Then $I_0,\dotsc,I_B$ are non-empty.
Pick vertices $v_0,\dotsc,v_B$ uniformly at random from $I_0,\dotsc,I_B$ respectively. Since, for every $i\neq j$,
the pair $v_iv_j$ is an edge with probability at least $1-\delta$, it follows that
$v_0,\dotsc,v_B$ is a clique with probability at least $1-\binom{B+1}{2}\delta>0$. In particular,
$G$ then contains a clique of size $B+1>1/\beta+1$, contrary to Lemma~\ref{lem:homog}. The contradiction
shows that $\abs{P}\leq BMd$, completing the proof of Theorem~\ref{thm:main}.

\section{Open problems}\label{sec:open}
\begin{itemize}
\item I know of only one asymptotic lower bound on $N_L$. It is a version of \cite[Proposition~5.12]{greaves_koolen_munemasa_szollosi}
that is also implicit in the bound for $N_{1/3}(d)$ in \cite{lemmens_seidel}. Denote by $I_n$ the identity matrix of size $n$, and by $J_n$ the all-one
matrix of size~$n$. Then the matrix $M=(r-1)I_{rt}-(J_r-I_r)\otimes I_t$ is a positive semidefinite matrix
of nullity $t$, it has $(r-1)$'s on the diagonal, and its off-diagonal entries are $0$ and $-1$.
Hence, $\tfrac{1}{r-1+\tau}(M+\tau J_{rt})$ is a Gram matrix of a $\{-\tfrac{1-\tau}{r-1+\tau},\tfrac{\tau}{r-1+\tau}\}$-code
in $\R^{(r-1)t+1}$ of size~$rt$. So, $N_L(d)\geq \frac{r}{r-1}d+O(1)$ for $L= \{-\tfrac{1-\tau}{r-1+\tau},\tfrac{\tau}{r-1+\tau}\}$.
For $\tau=1/2$, this yields a family of equiangular lines.
The results in \cite{lemmens_seidel,neumaier,greaves_koolen_munemasa_szollosi} suggest that this bound is sharp.
\begin{conjecture}
For an integer $r\geq 2$, the maximum number of equiangular lines with angle $\arccos \frac{1}{2r-1}$ is $N_{1/(2r-1)}(d)=\frac{r}{r-1}d+O(1)$ as $d$ tends to infinity.
\end{conjecture}

In contrast, one can show that the bound implicit in the proof of Theorem~\ref{thm:main} is $2^{O(1/\beta^2)}d$.

\item Informally, it is natural to think of Theorem~\ref{thm:main} as a juxtaposition of two trivial results from Lemma~\ref{lem:homog}: 
$N_{[-1,-\beta]}(d)= O(1)$ and $N_{\{\alpha\}}(d)=O(d)$. Since $N_{\{\alpha_1,\dotsc,\alpha_k\}}(d)=O(d^k)$ for
any real numbers $\alpha_1,\dotsc,\alpha_k$ (see~\cite[Proposition 1]{bukh_lmatrices}) this motivates the following conjecture.
\begin{conjecture}\label{conj:sudakov_keevash}
Suppose $\alpha_1,\dotsc,\alpha_k$ are any $k$ real numbers, and $L=[-1,-\beta]\cup\{\alpha_1,\dotsc,\alpha_k\}$. 
Then $N_L(d)\leq c_{\beta,k} d^k$.
\end{conjecture}
It is conceivable that in this case even $N_L(d)\leq c_\beta N_{\{\alpha_1,\dotsc,\alpha_k\}}(d)$ might be true.

\textbf{Added in revision:} Conjecture~\ref{conj:sudakov_keevash} has been resolved by Keevash and Sudakov \cite{keevash_sudakov}

\item I cannot rule out the possibility that \emph{for a fixed $\alpha$}
the size of any $[-1,0)\cup\{\alpha\}$-code is at most linear in the dimension.
\end{itemize}\smallskip

\textbf{Acknowledgments.} I am grateful to James Cummings, Hao Huang and Humberto Naves for inspirational discussions.
I am thankful to Joseph Briggs for careful reading, and for finding a mistake in an earlier version of this paper.
I also benefited from the constructive comments of two referees. 

\bibliographystyle{plain}
\bibliography{equiangular}

\end{document}